\documentclass{amsart}

\usepackage{amsfonts,amsmath,amssymb,mathrsfs,verbatim,cite}

\numberwithin{equation}{section}

\newtheorem{lemma}{Lemma}[section]
\newtheorem{theorem}{Theorem}[section]
\newtheorem{corollary}{Corollary}[section]
\newtheorem{proposition}{Proposition}[section]

\newcommand{\abs}[1]{\lvert#1\rvert}
\newcommand{\X}{\mathbb X}
\newcommand{\Y}{\mathbb Y}

\newcommand{\F}{\mathbb F}

\newcommand{\Z}{\mathbb Z}
\newcommand{\la}{\lambda}
\newcommand{\lae}{\lambda_{\textup{e}}}
\newcommand{\lao}{\lambda_{\textup{o}}}

\newcommand{\muo}{\mu_{\textup{o}}}

\newcommand{\qbin}[2]{\genfrac{[}{]}{0pt}{}{#1}{#2}}

\begin{document}

\title[Hall--Littlewood symmetric functions]
{Rogers--Szeg\"o polynomials and 
Hall--Littlewood symmetric functions}

\author{S. Ole Warnaar}\thanks{
Work supported by the Australian Research Council.}
\address{Department of Mathematics and Statistics,
The University of Melbourne, VIC 3010, Australia}
\email{warnaar@ms.unimelb.edu.au}


\begin{abstract}
We use Rogers--Szeg\"o polynomials to unify some well-known
identities for Hall--Littlewood symmetric functions due to Macdonald
and Ka\-wa\-na\-ka.
\end{abstract}

\maketitle

\section{Introduction and summary of results}

Three classical identities for Schur functions are
\cite{Schur73,Littlewood50,Macdonald95}
\begin{subequations}\label{s3}
\begin{equation}\label{Bn}
\sum_{\la} s_{\la}(x)=
\prod_{i\geq 1}\frac{1}{1-x_i}
\prod_{i<j}\frac{1}{1-x_ix_j}
\end{equation}
and
\begin{equation}\label{Cn}
\sum_{\substack{\la \\ \la\text{ even}}}
s_{\la}(x)=\prod_{i\geq 1}\frac{1}{1-x_i^2}
\prod_{i<j}\frac{1}{1-x_ix_j}
\end{equation}
and
\begin{equation}\label{Dn}
\sum_{\substack{\la \\ \la'\text{ even}}}
s_{\la}(x)=
\prod_{i<j}\frac{1}{1-x_ix_j}.
\end{equation}
\end{subequations}
Here $\la$ denotes a partition, $\la'$ its conjugate,
and the condition ``$\la$ even'' (or ``$\la'$ even'')
implies that all parts of $\la$ (or all parts of
$\la'$) must be even. Furthermore,
$s_{\la}(x)=s_{\la}(x_1,x_2,\dots)$ is a Schur function 
of a finite or infinite number of variables. 

When $x=(x_1,\dots,x_n)$ the identities \eqref{Bn}--\eqref{Dn}
may be viewed as reciprocals of Weyl denominator 
formulas; the latter expressing the products
\begin{equation*}
\prod_{i=1}^n (1-x_i)\prod_{1\leq i<j\leq n}(1-x_ix_j),
\qquad
\prod_{i=1}^n (1-x_i^2)\prod_{1\leq i<j\leq n}(1-x_ix_j)
\end{equation*}
and
\begin{equation*}
\prod_{1\leq i<j\leq n}(1-x_ix_j)
\end{equation*}
as sums over the B$_n$, C$_n$ and D$_n$ Weyl 
groups~\cite{Humphreys72}. 

Probably the most important application of \eqref{s3} was given
by Macdonald, who used the bounded form 
\begin{equation}\label{Bfin}
\sum_{\substack{\la \\ \la_1\leq k}}
s_{\la}(x_1,\dots,x_n)=\frac{\det\bigl(x_i^{j-1}-x_i^{2n+k-j}\bigr)}
{\prod_{i=1}^n (1-x_i)\prod_{1\leq i<j\leq n}(x_i-x_j)(1-x_ix_j)}
\end{equation}
of \eqref{Bn} to prove the famous MacMahon conjecture 
in the theory of plane partitions \cite{Bressoud99,Macdonald95}.

Formulae that incorporate all three Schur function identities
\eqref{s3} were recently found by Bressoud \cite{Bressoud00},
Ishikawa and Wakayama \cite{IW99},
and by Jouhet and Zeng \cite{JZ01}. 
If $m_i(\la)$ denotes the multiplicity of the part $i$ in
$\la$, i.e., $m_i(\la)=\la'_i-\la'_{i+1}$, then the
Bressoud--Ishikawa--Wakayama identity states that
\begin{equation}\label{IWid}
\sum_{\la} f_{\la}(a,b) s_{\la}(x)=
\prod_{i\geq 1}\frac{1}{(1-ax_i)(1-bx_i)}
\prod_{i<j}\frac{1}{1-x_ix_j},
\end{equation}
where
\begin{equation*}
f_{\la}(a,b)=
\prod_{j\text{ odd}}
\frac{a^{m_j(\la')+1}-b^{m_j(\la')+1}}{a-b}
\prod_{j\text{ even}}
\frac{1-(ab)^{m_j(\la')+1}}{1-ab}.
\end{equation*}
Similarly, the Jouhet--Zeng formula asserts that
\begin{equation}\label{JZid}
\sum_{\la} f_{\la'}(a,b) s_{\la}(x)=
\prod_{i\geq 1}\frac{(1+ax_i)(1+bx_i)}{(1-x_i)(1+x_i)}
\prod_{i<j}\frac{1}{1-x_ix_j}.
\end{equation}
For $b=0$ \eqref{IWid} and \eqref{JZid} reduce to identities
of Littlewood \cite{Littlewood50} combining \eqref{Bn}
and \eqref{Dn}, or \eqref{Bn} and \eqref{Cn}, respectively. 
Even more general formulae than \eqref{IWid} and \eqref{JZid}, 
which will not play a role in the present paper, may be found 
in \cite{JZ01}. A $\lambda$-ring approach to the
above results may be found in \cite{Lascoux03}.

\medskip

An important generalization of the Schur functions is given by
the Hall--Little\-wood symmetric functions $P_{\la}(x;t)$.
Here $t$ is an additional scalar variable
such that $P_{\la}(x;0)=s_{\la}(x)$.
Employing the Hall--Littlewood functions,
Macdonald \cite{Macdonald95} gave the following four
generalizations of the identities \eqref{Bn}--\eqref{Dn}:
\begin{subequations}\label{i16}
\begin{equation}\label{i1}
\sum_{\la} P_{\la}(x;t)=
\prod_{i\geq 1}\frac{1}{1-x_i}
\prod_{i<j}\frac{1-tx_ix_j}{1-x_ix_j}
\end{equation}
and
\begin{equation}\label{i2}
\sum_{\substack{\la \\ \la\text{ even}}} 
P_{\la}(x;t)=\prod_{i\geq 1}\frac{1}{1-x_i^2}
\prod_{i<j}\frac{1-tx_ix_j}{1-x_ix_j}
\end{equation}
and
\begin{equation}\label{i3}
\sum_{\substack{\la \\ \la'\text{ even}}} 
c_{\la}(t) P_{\la}(x;t)=
\prod_{i<j}\frac{1-tx_ix_j}{1-x_ix_j}
\end{equation}
and
\begin{equation}\label{i4}
\sum_{\la} d_{\la}(t)
P_{\la}(x;t)=
\prod_{i\geq 1}\frac{1-tx_i}{1-x_i}
\prod_{i<j}\frac{1-tx_ix_j}{1-x_ix_j},
\end{equation}
where for $\la'$ even (so that $m_i(\la)$ is even)
\begin{equation*}
c_{\la}(t)=\prod_{i\geq 1}(1-t)(1-t^3)\cdots (1-t^{m_i(\la)-1}),
\end{equation*}
and for general $\la$
\begin{equation*}
d_{\la}(t)=\prod_{i\geq 1}
(1-t)(1-t^3)\cdots (1-t^{2\lceil m_i(\la)/2\rceil-1})
\end{equation*}
with $\lfloor \cdot \rfloor$ and $\lceil \cdot \rceil$
the usual floor (or integer part) and ceiling functions.

Recently, Kawanaka \cite{Kawanaka99}
added two further identities to the list as follows.
For $\la$ a partition, let $\lae$ and $\lao$ be the partitions
containing the even parts and the odd parts of $\la$ respectively.
For example, if $\la=(4,3,3,2,1,1,1)$ then
$\lae=(4,2)$ and $\lao=(3,3,1,1,1)$. As usual
$l(\la)$ denotes the length of the partition $\la$ (that is,
the number of nonzero parts).
Then the Kawanaka identities correspond to the sums
\begin{equation}\label{i5}
\sum_{\la} e_{\la}(t) P_{\la}(x;t)=
\prod_{i\geq 1}\frac{1+t^{1/2}x_i}{1-x_i}
\prod_{i<j}\frac{1-tx_ix_j}{1-x_ix_j}
\end{equation}
and
\begin{equation}\label{i6}
\sum_{\substack{\la \\ (\lao)'\text{ even}}} 
f_{\la}(t) P_{\la}(x;t)=
\prod_{i\geq 1}\frac{1-tx_i^2}{1-x_i^2}
\prod_{i<j}\frac{1-tx_ix_j}{1-x_ix_j},
\end{equation}
\end{subequations}
where 
\begin{equation*}
e_{\la}(t)=\prod_{i\geq 1}(1+t^{1/2})(1+t)\cdots (1+t^{m_i(\la)/2}),
\end{equation*}
and, for $(\lao)'$ even (so that the odd parts of $\la$ have
even multiplicity),
\begin{equation*}
f_{\la}(t)=t^{l(\la_0)/2}d_{\la}(t).
\end{equation*}

Like their Schur function counterparts the above Hall--Littlewood
identities have interesting applications.
For example, Kawanaka's identities have an interpretation in
terms of the representation theory of the general linear group 
over finite fields \cite{Kawanaka91,Kawanaka99}. 
In particular, \eqref{i5} encodes the
fact that the symmetric space $\text{GL}_n(\F_{p^2})/\text{GL}_n(\F_p)$
(where $\text{GL}_n(\F_p)$ is the general linear group over a finite 
field of $p$ elements) is multiplicity free.
Similarly, \eqref{i6}, asserts that the symmetric space
$\text{GL}_{2n}(\F_p)/\text{Sp}_{2n}(\F_p)$
(with $\text{Sp}_{2n}$ the symplectic group)
is multiplicity free.

Another nice application follows by again considering the bounded 
versions of the identities of \eqref{i16}, see e.g., 
\cite{IJZ04,JZ04,Macdonald95,Stembridge90}.
For example, \eqref{i1} has the following bounded form generalizing
\eqref{Bfin}. Let $x=(x_1,\dots,x_n)$ and
\begin{equation*}
\Phi(x;t)=\prod_{i=1}^n\frac{1}{1-x_i}
\prod_{1\leq i<j\leq n} \frac{1-tx_ix_j}{1-x_ix_j}.
\end{equation*}
Then for $k$ a positive integer
\begin{equation*}
\sum_{\substack{\la \\ \la_1\leq k}} P_{\la}(x;t)
=\sum_{\epsilon\in\{\pm 1\}^n} \Phi(x^{\epsilon};t)
\prod_{i=1}^n x_i^{k(1-\epsilon_i)/2},
\end{equation*}
where $x^{\epsilon}=(x_1^{\epsilon_1},\dots,x_n^{\epsilon_n})$
\cite[pp. 232--234]{Macdonald95}.
Making the principal specialization $x=(z,zt,\dots,zt^{n-1})$
(and then replacing $t$ by $q$)
leads to interesting $q$-series identities.
The most important ones being the famous Rogers--Ramanujan identities
--- arising from the bounded form of \eqref{i2} due to Stembridge 
\cite{Stembridge90}.

\bigskip

Given the identities \eqref{i1}--\eqref{i6} and their striking
similarity, an obvious question is whether one can understand
all six as special cases of a master identity for Hall--Littlewood
functions. We will answer this question in the affirmative
in the form of Theorem~\ref{thmM} below, generalizing
the Jouhet--Zeng identity \eqref{JZid} to the level
of Hall--Littlewood functions.

For $m$ a nonnegative integer let $H_m(z;t)$ be the
Rogers--Szeg\"o polynomial \cite[Ch. 3, Examples 3--9]{Andrews76}
\begin{equation}\label{RS}
H_m(z;t)=\sum_{i=0}^m z^i \qbin{m}{i}_t.
\end{equation}
Here $\qbin{m}{n}_t$ is the usual $t$-binomial coefficient:
\begin{equation*}
\qbin{n}{m}_t=
\begin{cases}
\displaystyle
\frac{(t^{n-m+1};t)_m}{(t;t)_m} & \text{for $m\geq 0$,} \\[3mm]
0 & \text{otherwise,}
\end{cases}
\end{equation*}
where $(t;t)_0=1$ and $(t;t)_n=\prod_{i=1}^n (1-t^i)$ are
$t$-shifted factorials.
We extend the definition of the Rogers--Szeg\"o polynomials to
partitions $\la$ by
\begin{equation}\label{hla}
h_{\la}(z;t)=\prod_{i\geq 1} H_{m_i(\la)}(z;t).
\end{equation}
For example $h_{(3,2,2,1)}=H_1^2 H_2$.

\begin{theorem}\label{thmM}
The following formal identity holds:
\begin{multline}\label{master}
\sum_{\la} 
a^{l(\lao)} h_{\lae}(ab;t)h_{\lao}(b/a;t)
P_{\la}(x;t) \\
=\prod_{i\geq 1}\frac{(1+ax_i)(1+bx_i)}{(1-x_i)(1+x_i)}
\prod_{i<j}\frac{1-tx_ix_j}{1-x_ix_j}.
\end{multline}
\end{theorem}
It is important to note that the left-hand side satisfies
the necessary symmetry under interchange of $a$ and $b$.
{}From the definition of the Rogers--Szeg\"o polynomials
it readily follows that $H_m(z^{-1};t)=z^{-m}H_m(z;t)$.
Hence, since $\sum_i m_i(\la)=l(\la)$,
\begin{equation*}
h_{\la}(z^{-1};t)=z^{-l(\la)}h_{\la}(z;t).
\end{equation*}
Applying this to $h_{\lao}(b/a;t)$ in \eqref{master} shows
that the left is invariant under the interchange of $a$ and $b$.

When $t=1$ the Hall--Littlewood functions reduce to the
monomial symmetric functions, i.e., $P_{\la}(x;1)=m_{\la}(x)$.
Since $h_{\la}(z;1)=(1+z)^{l(\la)}$ this implies the elegant
summation
\begin{equation*}
\sum_{\la} (1+ab)^{l(\lae)}(a+b)^{l(\lao)} m_{\la}(x) 
=\prod_{i\geq 1}\frac{(1+ax_i)(1+bx_i)}{(1-x_i)(1+x_i)}.
\end{equation*}

In the following we will show how all six identities stated in \eqref{i16} 
follow from \eqref{master}. If we take $a=1$, use
$h_{\lae}(b;t)h_{\lao}(b;t)=h_{\la}(b;t)$ and finally replace
$b\to a$ we obtain our first corollary.
\begin{corollary}\label{cor1}
There holds
\begin{equation}\label{eqcor1}
\sum_{\la}
h_{\la}(a;t)P_{\la}(x;t) 
=\prod_{i\geq 1}\frac{1+ax_i}{1-x_i}
\prod_{i<j}\frac{1-tx_ix_j}{1-x_ix_j}.
\end{equation}
\end{corollary}
The following explicit evaluations for the Rogers--Szeg\"o polynomials
are known, see e.g., \cite{Andrews76,BW05}:
\begin{subequations}\label{spec}
\begin{align}
H_m(0;t)&=1 \\
H_m(-1;t)&=\begin{cases}(t;t^2)_{m/2} & \text{$m$ even} \\
0 & \text{$m$ odd} \end{cases} \\
H_m(-t;t)&=(t;t^2)_{\lceil m/2\rceil}=
(t;t)_m/(t^2;t^2)_{\lfloor m/2\rfloor}\label{spec3} \\[2mm]
H_m(t^{1/2};t)&=(-t^{1/2};t^{1/2})_m.
\end{align}
\end{subequations}
This immediately yields (in exactly that order)
\eqref{i1}, \eqref{i3}--\eqref{i5}.
We note in particular that by taking $a=-1$ in \eqref{eqcor1}
the summand vanishes unless all $m_i(\la)$ are even. That
is, all parts of $\la$ must have even multiplicity, or equivalently,
$\la'$ must be even.

Next we consider the case $b=-a$ of \eqref{master}. Using \eqref{spec3}
and making the replacement $a^2\to a$ this gives our second corollary.
\begin{corollary}
There holds
\begin{multline}\label{master2}
\sum_{\substack{\la \\ (\lao)' \textup{ even}}}
a^{l(\lao)/2} h_{\lae}(-a;t)
\biggl(\:\prod_{i\geq 1} (t;t^2)_{m_i(\lao)/2}\biggr) P_{\la}(x;t) \\
=\prod_{i\geq 1}\frac{1-ax_i^2}{1-x_i^2}
\prod_{i<j}\frac{1-tx_ix_j}{1-x_ix_j}.
\end{multline}
\end{corollary}
As remarked before, $(\lao)'$ being even is equivalent to
the odd parts of $\la$ having even multiplicity.
We also note that the product on the left-hand side
may be replaced by the equivalent
\begin{equation*}
\prod_{i\geq 1} (t;t^2)_{m_{2i-1}(\la)/2}.
\end{equation*}

Using three of the four specializations of \eqref{spec}
gives \eqref{i2}, \eqref{i3} and \eqref{i6}.
This shows that a seventh identity, corresponding
to \eqref{master2} with $a=-t^{1/2}$
has actually been missing from the literature:
\begin{equation*}
\sum_{\substack{\la \\ (\lao)' \text{ even}}} k_{\la}(t) P_{\la}(x;t)
=\prod_{i\geq 1}\frac{1+t^{1/2}x_i^2}{1-x_i^2}
\prod_{i<j} \frac{1-tx_ix_j}{1-x_ix_j},
\end{equation*}
with
\begin{equation*}
k_{\la}(t)=(-t^{1/2})^{l(\lao)/2}
\prod_{i\geq 1} (-t^{1/2};t^{1/2})_{m_i(\lae)}(t;t^2)_{m_i(\lao)/2}.
\end{equation*}

A further interesting special case of the theorem 
arises after taking $b=0$.
\begin{corollary}\label{cor3}
There holds
\begin{equation}\label{LW2}
\sum_{\la} a^{l(\lao)} P_{\la}(x;t)=
\prod_{i\geq 1}\frac{1+ax_i}{1-x_i^2}
\prod_{i<j}\frac{1-tx_ix_j}{1-x_ix_j}.
\end{equation}
\end{corollary}
In the Schur case this reduces to 
the Littlewood formula mentioned after \eqref{JZid},
combining \eqref{Bn} and \eqref{Cn}.

\bigskip

Observing that
\begin{equation*}
H_m(z;0)=1+z+\cdots+z^m=\frac{1-z^{m+1}}{1-z}
\end{equation*}
it readily follows that \eqref{master} simplifies to
\eqref{JZid} when $t=0$.
The reader may wonder whether there perhaps is
a companion to Theorem~\ref{thmM} extending
\eqref{IWid} in much the same way.
It is certainly possible (see \eqref{ms})
to obtain a formula of the form
\begin{equation}\label{fail}
\sum_{\la} C_{\la}(a,b;t) P_{\la}(x;t)=
\prod_{i\geq 1}\frac{(1-atx_i)(1-btx_i)}{(1-ax_i)(1-bx_i)}
\prod_{i<j}\frac{1-tx_ix_j}{1-x_ix_j}.
\end{equation} 
However, for general $a$ and $b$ the rational
function $C_{\la}(a,b;t)$ 
does not possess nice characteristics (like factorisation),
and we dismiss \eqref{fail} for being insufficiently interesting.
Only for $b=0$ we have a result elegant enough (although not very
deep) to be stated explicitly:
\begin{equation}\label{LW1}
\sum_{\la} a^{l((\la')_{\textup{o}})} d_{\la}(t)P_{\la}(x;t)=
\prod_{i\geq 1}\frac{1-a t x_i}{1-a x_i}
\prod_{i<j}\frac{1-tx_ix_j}{1-x_ix_j}.
\end{equation}
Here $(\la')_{\textup{o}}$ is the odd part of the conjugate of
$\la$ (so that $l((\la')_{\textup{o}})$ is the number of odd
columns of the diagram of $\la$) and $d_{\la}(t)=h_{\la}(-t;t)$
as before.
In the Schur case \eqref{LW1} reduces to 
the Littlewood formula mentioned after \eqref{JZid},
combining \eqref{Bn} and \eqref{Dn}.

It may perhaps seem surprising that at the level of Schur
functions a pair of equally elegant formulae exists but that
only one of these admits an appealing generalization to Hall--Littlewood
functions. The explanation for this is however easily given.
Let $\Lambda$ be the ring of symmetric functions and
$\omega:~\Lambda\to\Lambda$ the involution defined by
\begin{equation}\label{omega}
\omega(s_{\la})=s_{\la'}.
\end{equation}
\begin{lemma}\label{lemomega}
Applying $\omega$ to \eqref{IWid} yields \eqref{JZid}.
\end{lemma}
Hence \eqref{IWid} and \eqref{JZid} may really be viewed as
one and the same identity.
Since no ``good'' $t$-analogue of $\omega$ exists
there is no guarantee for Hall--Littlewood identities to
come in pairs also.

\bigskip
Finally we mention some further results related to \eqref{master}.
The first concerns the bounded form of Theorem~\ref{thmM},
or, to be precise, our failure to find this in full
generality.
At present we have only been able to find the bounded analogue of
\eqref{eqcor1} as follows.

For $x=(x_1,\dots,x_n)$ define
\begin{equation*}
\Phi(x;a,t)=\prod_{i=1}^n\frac{1+ax_i}{1-x_i}
\prod_{1\leq i<j\leq n} \frac{1-tx_ix_j}{1-x_ix_j}.
\end{equation*}
For $k$ a positive integer also define a bounded version of
$h_{\la}(z;t)$ by
\begin{equation*}
h_{\la;k}(z;t)=
\prod_{i=1}^{k-1} H_{m_i(\la)}(z;t).
\end{equation*}

For example $h_{(3,2,2,1);1}=1$,
$h_{(3,2,2,1);2}=H_1$,
$h_{(3,2,2,1);3}=H_1H_2$ and
$h_{(3,2,2,1);k}=H_1^2 H_2$ for $k\geq 4$.

\begin{proposition}\label{prop1}
For $k$ a positive integer and $x=(x_1,\dots,x_n)$
there holds
\begin{equation}\label{eqprop}
\sum_{\substack{\la \\ \la_1\leq k}}
h_{\la;k}(z;t) P_{\la}(x;t)
=\sum_{\epsilon\in\{\pm 1\}^n} \Phi(x^{\epsilon};a,t)
\prod_{i=1}^n x_i^{k(1-\epsilon_i)/2}.
\end{equation}
\end{proposition}
For $k=1$ this can be simplified since
\begin{equation*}
\sum_{\substack{\la \\ \la_1\leq 1}} P_{\la}(x;t)
=\sum_{r=0}^{\infty}P_{(1^r)}(x;t)
=\sum_{r=0}^{\infty}e_r(x)
=\prod_{i\geq 1}(1+x_i),
\end{equation*}
with $e_r$ the $r$th elementary symmetric function.
Hence
\begin{equation*}
\sum_{\epsilon\in\{\pm 1\}^n} \Phi(x^{\epsilon};a,t)
\prod_{i=1}^n x_i^{(1-\epsilon_i)/2}=\prod_{i=1}^n(1+x_i).
\end{equation*}

Assuming the specialization $x=(z,zt,\dots,zt^{n-1})$,
replacing $t$ by $q$, and letting $n$ tend to infinity
yields the $b=0$ case of the next proposition.
\begin{proposition}
There holds
\begin{align*}
\sum_{\substack{\la \\ \la_1\leq k}} &
z^{\abs{\la}} (b;q^{-1})_{\la_1'}
h_{\la;k}(a;q) P_{\la}(1,q,q^2,\dots;q) \\
&=\frac{(bz^2,-z,-az;q)_{\infty}}{(z^2,-bz,-abz;q)_{\infty}} \\
& \qquad \qquad \times
\sum_{r=0}^{\infty} 
(-1)^r a^r z^{kr}q^{(k+1)\binom{r}{2}}\,
\frac{1-z^2q^{2r-1}}{1-z^2q^{-1}}\;
\frac{(b;q^{-1})_r (z^2/q,-z/a;q)_r}{(q,-az,bz^2;q)_r}.
\end{align*}
\end{proposition}
Because it lies somewhat outside the scope of the present paper
we will not prove this $q$-series identity here. (For $b=0$ it
of course follows from \eqref{eqprop}.) 
As one application let us take $b=0$ and assume that $z=q^{1/2}$
or $z=q$ but $-az\neq q$. Some simple manipulations then give
\begin{multline*}
\sum_{\substack{\la \\ \la_1\leq k}}
z^{\abs{\la}} h_{\la;k}(a;q) P_{\la}(1,q,q^2,\dots;q) \\
=\frac{(-z,-az;q)_{\infty}}{(q;q)_{\infty}}
\sum_{r=-\infty}^{\infty} (-1)^r a^r z^{kr} q^{(k+1)\binom{r}{2}}
\frac{(-z/a;q)_r}{(-az;q)_r}.
\end{multline*}
When $z=q$ and $a=1$ the right may be written as a product 
using Jacobi's triple product identity, so that
we find the Rogers--Ramanujan-type identity
\begin{equation*}
\sum_{\substack{\la \\ \la_1\leq k}}
q^{\abs{\la}} h_{\la;k}(1;q) P_{\la}(1,q,q^2,\dots;q)
=\frac{(-q;q)_{\infty}^2(q,q^k,q^{k+1};q^{k+1})_{\infty}}{(q;q)_{\infty}}.
\end{equation*}

\bigskip

Finally we mention that for general $a$ and $b$ the generalization
of Theorem~\ref{thmM} to Macdonald polynomials lacks the necessary
elegance, and only \eqref{LW2} and \eqref{LW1} admit simple
$q$-generalizations.

Let $P_{\la}(x;q,t)$ be Macdonald's symmetric function
and let $b_{\la}^{\textup{oa}}(q,t)$ and
$b_{\la}^{\textup{el}}(q,t)$ be the rational functions defined in
\eqref{oael} of the next section.
\begin{proposition}\label{prop3}
The following formal identities hold:
\begin{equation}\label{Pb}
\sum_{\la} a^{l(\lao)} b_{\la}^{\textup{oa}}(q,t) P_{\la}(x;q,t)=
\prod_{i\geq 1} \frac{(1+ax_i)(q t x_i^2;q^2)_{\infty}}
{(x_i^2;q^2)_{\infty}}
\prod_{i<j}\frac{(tx_ix_j;q)_{\infty}}{(x_ix_j;q)_{\infty}}
\end{equation}
and
\begin{equation}\label{Pc}
\sum_{\la} a^{l((\la')_{\textup{o}})} 
b_{\la}^{\textup{el}}(q,t) P_{\la}(x;q,t)=
\prod_{i\geq 1}\frac{(a t x_i;q)_{\infty}}{(a x_i;q)_{\infty}}
\prod_{i<j}\frac{(tx_ix_j;q)_{\infty}}{(x_ix_j;q)_{\infty}}.
\end{equation}
\end{proposition}

\bigskip

In the next section we give a brief introduction to 
Hall--Littlewood functions,
Section~\ref{sec3} contains a proof of the claims of the first section,
and, finally, in Section~\ref{sec4} we restate some of our results
in the language of $\lambda$-rings.

\section{Hall--Littlewood functions}\label{sec2}
Let $\lambda=(\lambda_1,\lambda_2,\dots)$ be a partition, 
i.e., $\lambda_1\geq \lambda_2\geq \dots$ with finitely many 
$\lambda_i$ unequal to zero.
The length and weight of $\lambda$, denoted by 
$l(\lambda)$ and $\abs{\lambda}$, are the number and sum
of the non-zero $\lambda_i$ (called parts), respectively.
The unique partition of weight zero is denoted by $0$, and
the multiplicity of the part $i$ in the partition $\la$ is denoted
by $m_i(\lambda)$.

We identify a partition with its (Young) diagram or Ferrers graph in the
usual way, and the conjugate $\lambda'$ of $\lambda$ is the partition obtained by reflecting the diagram of $\lambda$ in the main diagonal.
Hence $m_i(\lambda)=\la_i'-\la_{i+1}'$.

If $\la$ and $\mu$ are two partitions then $\mu\subset\la$ iff
$\la_i\geq \mu_i$ for all $i\geq 1$, i.e., the diagram of $\la$ 
contains the diagram of $\mu$. If $\mu\subset\la$ then 
the skew-diagram $\la-\mu$ denotes the set-theoretic difference 
between $\la$ and $\mu$, and $\abs{\la-\mu}=\abs{\la}-\abs{\mu}$.
A skew diagram $\theta$ is a horizontal/vertical $r$-strip if 
it contains exactly $r$ squares, i.e., $\abs{\theta}=r$, and has at
most one square in each of its columns/rows.
For example, if $\la=(6,3,3,1)$ and $\mu=(4,3,1)$ then
$\la-\mu$ is a horizontal $5$-strip and $\la'-\mu'$ a vertical $5$-strip.

Let $s=(i,j)$ be a square in the diagram of $\lambda$. Then
$a(s)$, $a'(s)$, $l(s)$ and $l'(s)$ are the arm-length, arm-colength,
leg-length and leg-colength of $s$, defined by
\begin{subequations}\label{aapllp}
\begin{align}
a(s)&=\lambda_i-j,  & a'(s)&=j-1 \\
l(s)&=\lambda'_j-i,  & l'(s)&=i-1.
\end{align}
\end{subequations}
{}From this we may define several standard
rational functions on partitions:
\begin{equation*}
b_{\la}(s;q,t)=\frac{1-q^{a(s)}t^{l(s)+1}}{1-q^{a(s)+1}t^{l(s)}},
\end{equation*}
and
\begin{align}\label{oael}
b_{\la}^{\text{el}}(q,t)&=
\prod_{\substack{s\in\la \\ l(s) \text{ even}}} b_{\la}(s;q,t), &
b_{\la}^{\text{oa}}(q,t)&=
\prod_{\substack{s\in\la \\ a(s) \text{ odd}}} b_{\la}(s;q,t), 
\\[2mm]
b_{\la}^{\text{ol}}(q,t)&=
\prod_{\substack{s\in\la \\ l(s) \text{ odd}}} b_{\la}(s;q,t), &
b_{\la}^{\text{ea}}(q,t)&=
\prod_{\substack{s\in\la \\ a(s) \text{ even}}} b_{\la}(s;q,t).
\notag
\end{align}
Obviously,
\begin{equation}\label{b1}
b_{\la}(q,t):=b_{\la}^{\text{el}}(q,t)b_{\la}^{\text{ol}}(q,t)=
b_{\la}^{\text{ea}}(q,t)b_{\la}^{\text{oa}}(q,t)
\end{equation}
and
\begin{equation}\label{b2}
b_{\la'}^{\text{el}}(q,t)b_{\la}^{\text{ea}}(t,q)=
b_{\la'}^{\text{ol}}(q,t)b_{\la}^{\text{oa}}(t,q)=1.
\end{equation}

\medskip

Let $S_n$ be the symmetric group, $\Lambda_n=\Z[x_1,\dots,x_n]^{S_n}$ 
be the ring of symmetric polynomials in $n$ 
independent variables and $\Lambda$ the ring of symmetric functions
in countably many variables.

For $x=(x_1,\dots,x_n)$ and $\lambda$ a partition such that $l(\la)\leq n$
the Hall--Littlewood polynomial $P_{\lambda}(x;t)$ is defined by
\begin{equation}\label{Pdef} 
P_{\lambda}(x;t)
=\sum_{w\in S_n/S_n^{\la}} 
w\Bigl(x^{\la} \prod_{\la_i>\la_j}
\frac{x_i-tx_j}{x_i-x_j} \Bigr).
\end{equation}
Here $S_n^{\la}$ is the subgroup of $S_n$ consisting of those
permutations that leave $\la$ invariant, and $w(f(x))=f(w(x))$.
When $l(\la)>n$, 
\begin{equation}\label{Pnul}
P_{\lambda}(x;t)=0.
\end{equation}

The Hall--Littlewood polynomials are symmetric polynomials in $x$,
homogeneous of degree $\abs{\la}$, with coefficients
in $\Z[t]$, and form a $\Z[t]$ basis of $\Lambda_n[t]$.
Thanks to the stability property $P_{\la}(x_1,\dots,x_n,0;t)=
P_{\la}(x_1,\dots,x_n;t)$ the Hall--Littlewood polynomials may be 
extended to the Hall--Littlewood functions in an infinite number of
variables $x_1,x_2,\dots$ in the usual way, to form a $\Z[t]$ basis of
$\Lambda[t]$.
The parameter $t$ in the Hall--Littlewood symmetric functions
serves to interpolate between the
Schur functions and monomial symmetric functions;
$P_{\la}(x;0)=s_{\la}(x)$ and $P_{\la}(x;1)=m_{\la}(x)$.
We also introduce a second Hall--Littlewood function
$Q_{\la}$ by
\begin{equation}\label{Qla}
Q_{\la}(x;t)=b_{\la}(t) P_{\la}(x;t),
\end{equation}
where $b_{\la}(t)=b_{\la}(0,t)=\prod_{i\geq 1}(t;t)_{m_i(\la)}$.
Then the Cauchy identity for Hall--Littlewood functions takes the form
\begin{equation}\label{Cauchy}
\sum_{\la} P_{\la}(x;t)Q_{\la}(y;t)=\prod_{i,j\geq 1}
\frac{1-tx_iy_j}{1-x_iy_j}.
\end{equation}

When $\la=(1^r)$ and $\la=(r)$ the Hall--Littlewood polynomials 
reduce to the $r$th elementary and $r$th complete symmetric 
functions
\begin{equation}\label{eh}
P_{(1^r)}=e_r \quad \text{ and }\quad P_{(r)}=h_r.
\end{equation}
These functions may be defined by their generating functions as
\begin{equation}\label{er}
\sum_{r=0}^{\infty}z^r e_r(x)=\prod_{i\geq 1}(1+tx_i)
\end{equation}
and
\begin{equation}\label{hr}
\sum_{r=0}^{\infty}z^r h_r(x)=\prod_{i\geq 1}\frac{1}{1-tx_i}.
\end{equation}
Since $e_r=s_{(1^r)}$ and $h_r=s_{(r)}$ we have
\begin{equation}
\omega(e_r)=h_r,
\end{equation} 
with $\omega$ the involution \eqref{omega}.

The Pieri formula for Hall--Littlewood polynomials
states that
\begin{equation}\label{Pieri}
P_{\mu}(x;t)e_r(x)=\sum_{\la}f_{\mu(1^r)}^{\la}(t)P_{\la}(x;t),
\end{equation}
where the coefficient $f_{\mu(1^r)}^{\la}(t)$
is zero unless $\mu\subset\la$ such that the skew 
diagram $\la-\mu$ is a horizontal $r$-strip.
An explicit expression for $f_{\mu(1^r)}^{\la}(t)$ is given by
\cite[p. 215]{Macdonald95}
\begin{equation}\label{Pieri2}
f_{\mu(1^r)}^{\la}(t)=
\prod_{i\geq 1} \qbin{\la_i'-\la_{i+1}'}{\la_i'-\mu_i'}_t
\qquad\text{for $\abs{\la-\mu}=r$}
\end{equation}
and zero otherwise.

The more general structure constants of the Hall--Littlewood functions
are defined by
\begin{equation}\label{struc}
P_{\mu}(x;t)P_{\nu}(x;t)=\sum_{\la}f_{\mu\nu}^{\la}(t)P_{\la}(x;t).
\end{equation}
These may be utilized to define the skew function $Q_{\la/\mu}$ by
\begin{equation}\label{skew}
Q_{\la/\mu}(x;t)=\sum_{\nu}f_{\mu\nu}^{\la}(t)Q_{\nu}(x;t).
\end{equation}

\section{Proofs}\label{sec3}

\subsection{Proof of Theorem~\ref{thmM}}

We first prove the $b=0$ case of the theorem, corresponding
to Corollary~\ref{cor3}, and then use this to obtain the theorem
for general $a$ and $b$.

Our point of departure is \eqref{i2}. 
Replacing the summation index $\la$ by $\mu$
and multiplying both sides by $\prod_i (1+ax_i)$ yields
\begin{equation}\label{L1}
\sum_{\substack{\mu \\ \mu\text{ even}}}
P_{\mu}(x;t)\prod_{i\geq 1}(1+ax_i)
=\prod_{i\geq 1}\frac{1+ax_i}{1-x_i^2}
\prod_{i<j}\frac{1-tx_ix_j}{1-x_ix_j}.
\end{equation}
By \eqref{er} we can expand the left-hand side of \eqref{L1} as
\begin{equation*}
\text{LHS}\eqref{L1}=
\sum_{r=0}^{\infty}
\sum_{\substack{\mu \\ \mu\text{ even}}}
a^r P_{\mu}(x;t)e_r(x).
\end{equation*}
Next we use the Pieri formula \eqref{Pieri} to rewrite this as
\begin{equation*}
\text{LHS}\eqref{L1}=
\sum_{r=0}^{\infty}
\sum_{\substack{\la,\mu \\ \mu\text{ even}}}
a^r f_{\mu(1^r)}^{\la}(t) P_{\la}(x;t).
\end{equation*}
Since $f_{\mu(1^r)}^{\la}(t)=0$ when $\abs{\la-\mu}\neq r$
this may also be written as
\begin{equation*}
\text{LHS}\eqref{L1}=
\sum_{\substack{\la,\mu \\ \mu\text{ even}}}
a^{\abs{\la-\mu}} f_{\mu(1^{\abs{\la-\mu}})}^{\la}(t) P_{\la}(x;t).
\end{equation*}
Since $f_{\mu(1^{\abs{\la-\mu}})}^{\la}(t)$ is zero unless
$\la-\mu$ is a vertical strip, only those partitions $\mu$
contribute to the sum for which $0\leq \la_i-\mu_i\leq 1$.
Combined with the fact that $\mu$ must be even this completely
fixes $\mu$ as $\mu_i=2\lfloor \la_i/2\rfloor$ (so that
$\abs{\la-\mu}=l(\lao)$, the number of parts of $\la$ of odd length).
For example if $\la=(7,5,5,4,3,1)$ then the
only contributing $\mu$ to the sum is $\mu=(6,4,4,4,2)$.
In terms of conjugate partitions this implies that
if $\la_i'>\la_{i+1}'$ then $\mu_i'=\la_{i+1}'$.
For the partitions in our example $\la'=(6,5,5,4,3,1,1)$ and
$\mu'=(5,5,4,4,1,1)$ and $\la_1>\la_2$ so that $\la_2=\mu_2$,
$\la_3>\la_4$ so that $\la_4=\mu_4$, et cetera. From \eqref{Pieri2} 
we infer that
\begin{equation}\label{f}
f_{\mu(1^{\abs{\la-\mu}})}^{\la}(t)=
\prod_{i\geq 1} \qbin{\la_i'-\la_{i+1}'}{\la_i'-\mu_i'}_t.
\end{equation}
By the above considerations regarding $\la$ and $\mu$,
we find that whenever an upper index of a $t$-binomial coefficient
in the above product is positive
the lower index must be zero. Hence we simplify to
\begin{equation*}
\text{LHS}\eqref{L1}=
\sum_{\la} a^{l(\lao)} P_{\la}(x;t).
\end{equation*}
Equating this with the right-hand side of \eqref{L1} completes
the proof of Corollary~\ref{cor3}.

Next we use \eqref{LW2} to prove the full theorem.
To this end we multiply both sides of \eqref{LW2} 
by $\prod_i (1+bx_i)$ and replace $\la$ by $\mu$ to get
\begin{equation}\label{L2}
\sum_{\mu} a^{l(\muo)} P_{\mu}(x;t)\prod_{i\geq 1}(1+bx_i)=
\prod_{i\geq 1}\frac{(1+ax_i)(1+bx_i)}{(1-x_i)(1+x_i)}
\prod_{i<j}\frac{1-tx_ix_j}{1-x_ix_j}.
\end{equation}
Following exactly the same steps as before,
again using \eqref{er}, \eqref{Pieri} and \eqref{f}, the
left-hand side may be rewritten as
\begin{equation*}
\text{LHS}\eqref{L2}=
\sum_{\la,\mu} a^{l(\muo)} b^{\abs{\la-\mu}} P_{\la}(x;t)
\prod_{i\geq 1} \qbin{\la_i'-\la_{i+1}'}{\la_i'-\mu_i'}_t.
\end{equation*}
Next we replace the sum over the partition $\mu$ 
by a sum over a sequence $k=(k_1,k_2,\dots)$ of nonnegative
integers as follows: $\mu_i'=\la_i'-k_i$.
Using $\la_i'-\la_{i+1}'=m_i(\la)$ and
\begin{align*}
l(\muo)&=\sum_{i\geq 1} m_{2i-1}(\mu) \\
&=\sum_{i\geq 1}(\mu_{2i-1}'-\mu_{2i}') \\
&=\sum_{i\geq 1}(\la_{2i-1}'-\la_{2i}'-k_{2i-1}+k_{2i}) \\
&=l(\lao)-\sum_{i\geq 1}(k_{2i-1}-k_{2i}) \\
&=l(\lao)+\sum_{i\geq 1}(-1)^i k_i,
\end{align*}
we then obtain
\begin{align*}
\text{LHS}\eqref{L2}&=
\sum_{\la,k} a^{l(\la)} P_{\la}(x;t)
\prod_{i\geq 1} a^{(-1)^i k_i} b^{k_i} \qbin{m_i(\la)}{k_i}_t \\
&=\sum_{\la} a^{l(\la)} P_{\la}(x;t)
\prod_{i\geq 1}\sum_{k_i=0}^{m_i(\la)}
a^{(-1)^i k_i} b^{k_i} \qbin{m_i(\la)}{k_i}_t.
\end{align*}
Factoring the product over $i$ into a product over even values of $i$
and a product over odd values of $i$ and then using that
\begin{equation*}
m_i(\lae)=\begin{cases} m_i(\la) &\text{if $i$ is even} \\
0 &\text{if $i$ is odd}
\end{cases} \quad\text{and}\quad
m_i(\lao)=\begin{cases} m_i(\la) &\text{if $i$ is odd} \\
0 &\text{if $i$ is even,}
\end{cases}
\end{equation*}
we obtain the further rewriting
\begin{equation*}
\text{LHS}\eqref{L2}
=\sum_{\la} a^{l(\lao)} P_{\la}(x;t)  
\prod_{i\geq 1} 
\sum_{k=0}^{m_i(\lae)} (ab)^k \qbin{m_i(\lae)}{k}_t
\sum_{k=0}^{m_i(\lao)} (b/a)^k \qbin{m_i(\lao)}{k}_t .
\end{equation*}
Finally, by \eqref{RS} and \eqref{hla}, this becomes
\begin{align*}
\text{RHS}\eqref{L2}
&=\sum_{\la} a^{l(\lao)} P_{\la}(x;t)  
\prod_{i\geq 1} H_{m_i(\lae)}(ab;t) H_{m_i(\lao)}(b/a;t) \\
&=\sum_{\la} a^{l(\lao)} 
h_{\lae}(ab;t) h_{\lao}(b/a;t) P_{\la}(x;t),
\end{align*}
completing the proof.

\subsection{Proof of Lemma~\ref{lemomega}}
Acting with $\omega$ on the left-hand side of \eqref{JZid}
yields
\begin{align*}
\omega(\text{LHS}\eqref{JZid})&=
\sum_{\la} f_{\la'}(a,b) \omega(s_{\la})(x) \\
&=\sum_{\la} f_{\la'}(a,b) s_{\la'}(x) \\
&=\sum_{\la} f_{\la}(a,b) s_{\la}(x) \\
&=\text{LHS}\eqref{IWid}).
\end{align*}
where in the second-last step we have changed the summation index from
$\la$ to its conjugate and used the fact
that summing over $\la$ is equivalent to summing over $\la'$.

Dealing with the right-hand side requires a few more steps
but is equally elementary. By \eqref{Cn} and \eqref{er} we have
\begin{equation*}
\text{RHS}\eqref{JZid})=
\sum_{u,v} \sum_{\substack{\la \\ \la\text{ even}}}
a^u b^v e_u(x)e_v(x)s_{\la}(x).
\end{equation*}
Therefore
\begin{align*}
\omega(\text{RHS}\eqref{JZid})&=
\sum_{u,v} \sum_{\substack{\la \\ \la\text{ even}}}
a^u b^v \omega(e_u)\omega(e_v)\omega(s_{\la})(x) \\
&=\sum_{u,v} \sum_{\substack{\la \\ \la\text{ even}}}
a^u b^v h_u(x)h_v(x)s_{\la'}(x) \\
&=\sum_{u,v} \sum_{\substack{\la \\ \la'\text{ even}}}
a^u b^v h_u(x)h_v(x)s_{\la}(x) \\
&=\sum_{u,v} 
a^u b^v h_u(x)h_v(x) \prod_{i<j}\frac{1}{1-x_ix_j},
\end{align*}
where the last equality follows from \eqref{Dn}.
Finally using \eqref{hr} we get
\begin{equation*}
\omega(\text{RHS}\eqref{JZid})=
\text{RHS}\eqref{IWid}.
\end{equation*}

\subsection{Proof of Proposition~\ref{prop1}}
The proof follows \cite{IJZ04,JZ04,Macdonald95,Stembridge90}
mutatis mutandis.

\subsection{Proof of Proposition~\ref{prop3}}
We will assume the reader is familiar with the theory of 
Macdonald polynomials. All notations and definitions used
in the proof may be found in Chapter VI of \cite{Macdonald95}.
Whenever possible we have indicated the precise page
in \cite{Macdonald95} where a particular result or definition may be
found.

\medskip

Both \eqref{Pb} and \eqref{Pc} may simply be proved using their
$a=0$ specializations established in \cite{Macdonald95}.
It is however more instructive to only prove \eqref{Pc} in this way,
and to obtain \eqref{Pb} by acting on the former with the
automorphism $\omega_{q,t}$ of $\Lambda_F$. This automorphism
acts on the Macdonald polynomials as \cite[p. 327]{Macdonald95}
\begin{equation*}
\omega_{q,t}P_{\la}(x;q,t)=Q_{\la'}(x;t,q),
\end{equation*}
where $Q_{\la}(x;q,t)=b_{\la}(q,t)P_{\la}(x;q,t)$.

\begin{proof}[Proof of \eqref{Pc}]
We may assume the $a=0$ case of \eqref{Pc} given by 
\cite[p. 349]{Macdonald95}
\begin{equation}\label{p349}
\sum_{\substack{\la \\ \la'\text{ even}}}
b_{\la}^{\textup{el}}(q,t) P_{\la}(x;q,t)=
\prod_{i<j}\frac{(tx_ix_j;q)_{\infty}}{(x_ix_j;q)_{\infty}}.
\end{equation}
Since \cite[p. 311]{Macdonald95}
\begin{equation}\label{gr}
\sum_{r=0}^{\infty} g_r(x;q,t) a^r =
\prod_{i\geq 1}\frac{(a t x_i;q)_{\infty}}{(a x_i;q)_{\infty}}
\end{equation}
this implies that
\begin{equation*}
\text{RHS}\eqref{Pc}=
\sum_{\substack{\mu,r \\ \mu'\text{ even}}}
b_{\mu}^{\text{el}}(q,t) a^r P_{\mu}(x;q,t)g_r(x;q,t).
\end{equation*}
By the Pieri formula \cite[p. 340]{Macdonald95}
\begin{equation*}
P_{\mu}(x;q,t) g_r(x;q,t)=
\sum_{\substack{\la \\ \la-\mu\text{ hor.  $r$-strip}}}
\varphi_{\la/\mu}(q,t) P_{\la}(x;q,t)
\end{equation*}
this becomes
\begin{equation*}
\text{RHS}\eqref{Pc}=
\sum_{\substack{\la,\mu \\ \mu' \text{ even} \\ \la-\mu \text{ hor. strip}}} 
a^{\abs{\la-\mu}} b_{\mu}^{\text{el}}(q,t)
\varphi_{\la/\mu}(q,t) P_{\la}(x;q,t).
\end{equation*}
Reasoning as before (see the proof of Theorem~\ref{thmM})
it follows that for given $\la$ the partition $\mu$ is uniquely fixed
as $\mu_i'=2\lfloor \la_i'/2\rfloor$.
Assuming such $\mu$ we thus obtain
\begin{equation*}
\text{RHS}\eqref{Pc}=\sum_{\la} 
a^{l((\la')_{\text{o}})} b_{\mu}^{\text{el}}(q,t)
\varphi_{\la/\mu}(q,t) P_{\la}(x;q,t).
\end{equation*}
Since \cite[p. 351]{Macdonald95}
\begin{equation*}
b_{\mu}^{\text{el}}(q,t)
\varphi_{\la/\mu}(q,t)=b_{\la}^{\textup{el}}(q,t)
\end{equation*}
(for $\mu_i'=2\lfloor \la_i'/2\rfloor$)
we arrive at 
\begin{equation*}
\text{RHS}\eqref{Pc}=\sum_{\la} 
a^{l((\la')_{\text{o}})} b_{\la}^{\text{el}}(q,t) P_{\la}(x;q,t)
\end{equation*}
completing the proof.
\end{proof}
A slightly different proof in the context of $\lambda$-rings will
be presented in the next section.

\begin{proof}[Proof of \eqref{Pb}]
Acting with $\omega_{q,t}$ on the left of \eqref{Pc} yields
\begin{align*}
\omega_{q,t}(\text{LHS}\eqref{Pc})
&= \sum_{\la} a^{l((\la')_{\text{o}})}
b_{\la}^{\text{el}}(q,t) Q_{\la'}(x;t,q) \\
&= \sum_{\la} a^{l((\la')_{\text{o}})}
b_{\la'}^{\text{el}}(q,t) b_{\la}(t,q)P_{\la}(x;t,q) \\
&= \sum_{\la} a^{l(\la_{\text{o}})}
b_{\la}^{\text{oa}}(t,q) P_{\la}(x;t,q),
\end{align*}
where the last equality follows by \eqref{b1} and \eqref{b2}.

On the other hand, by \eqref{gr} the right-hand side of \eqref{Pc}
may be written as
\begin{equation*}
\text{RHS}\eqref{Pc}=
\sum_{r=0}^{\infty} g_r(x;q,t) a^r \:
\prod_{i<j}\frac{(tx_ix_j;q)_{\infty}}{(x_ix_j;q)_{\infty}}.
\end{equation*}
Applying $\omega_{q,t}$ and
using \cite[p. 312]{Macdonald95}
\begin{equation*}
\omega_{q,t}(g_r(x;q,t))=e_r(x)
\end{equation*}
and \cite[p. 351]{Macdonald95}
\begin{equation*}
\omega_{q,t}\biggl(\:
\prod_{i<j}\frac{(tx_ix_j;q)_{\infty}}{(x_ix_j;q)_{\infty}}
\biggr)=
\prod_{i\geq 1} \frac{(q t x_i^2;t^2)_{\infty}}{(x_i^2;t^2)_{\infty}}
\prod_{i<j}\frac{(qx_ix_j;t)_{\infty}}{(x_ix_j;t)_{\infty}}
\end{equation*}
gives
\begin{align*}
\omega_{q,t}(\text{RHS}\eqref{Pc})&=
\sum_{r=0}^{\infty} e_r(x) a^r \:
\prod_{i\geq 1} \frac{(q t x_i^2;t^2)_{\infty}}{(x_i^2;t^2)_{\infty}}
\prod_{i<j}\frac{(qx_ix_j;t)_{\infty}}{(x_ix_j;t)_{\infty}} \\
&=\prod_{i\geq 1} \frac{(1+ax_i)(q t x_i^2;t^2)_{\infty}}
{(x_i^2;t^2)_{\infty}}
\prod_{i<j}\frac{(qx_ix_j;t)_{\infty}}{(x_ix_j;t)_{\infty}}.
\qedhere
\end{align*}
\end{proof}

\section{$\la$-rings}\label{sec4}
Lascoux recently revisited the Schur function identities of the 
introduction from the point of view of $\lambda$-rings \cite{Lascoux03}.
In this section we adopt Lascoux's approach, and restate some 
of our results in $\la$-ring (or plethystic) notation. 
For an introduction to symmetric functions and $\la$-rings 
we refer to \cite{Knutson73,Lascoux03b}.

Given two alphabets $\X$ and $\Y$ we denote by 
$\X+\Y$ and $\X\Y$ their disjoint union and Cartesian product.
Decomposing an alphabet as the sum of its letters,
we follow the convention of writing
$\X=\sum_{x\in\X}x$ instead of $\X=\sum_{x\in\X}\{x\}$.

The complete symmetric function
$h_r[\X-\Y]$ is defined by its generating series
\begin{equation}\label{hrla}
\sigma_z[\X-\Y]:=\frac{\prod_{y\in\Y} (1-zy)}{\prod_{x\in\X}(1-zx)}
=\sum_{r=0}^{\infty}z^r h_r[\X-\Y].
\end{equation}
Here we use the plethystic brackets to distinguish from our earlier notation
of \eqref{hr}. In particular, $h_r(x_1,x_2,\dots)=h_r[\X]$ and
$e_r(x_1,x_2,\dots)=e_r[\X]=(-1)^r h_r[-\X]$ for $\X=\{x_1,x_2,\dots\}$.
We also define $h_r[(1-q)\X/(1-t)]$ by 
\begin{equation*}
\prod_{x\in\X}\frac{(tzx;q)_{\infty}}{(zx;q)_{\infty}}=
\sum_{r\geq 0} z^r h_r[(1-t)\X/(1-q)].
\end{equation*}
(We mostly use this with $t=0$ and $q$ replaced by $t$.)
Then, by \eqref{hrla}, $\X/(1-t)=\X\{1,t,t^2,\dots\}$,
so that by Euler's $q$-exponential sum \cite[Equation (II.1)]{GR04}
\begin{equation}\label{Euler}
h_r[a/(1-t)]=\frac{a^r}{(t;t)_r}.
\end{equation}

For our present purposes it is important to
note that the Rogers--Szeg\"o polynomials actually 
arise as complete symmetric functions \cite[Exercise 2.22]{Lascoux03}:
\begin{equation*}
a^r H_r(b/a;t)=
\frac{h_r[\X/(1-t)]}{h_r[1/(1-t)]}
=(t;t)_r h_r[\X/(1-t)],\qquad \X=\{a,b\}.
\end{equation*}
Indeed, since $\X=\{a,b\}=a+b$, the factorization of 
the left-hand side of \eqref{hrla} implies the convolution
\begin{align}\label{RShr}
h_r[\X/(1-t)]&=\sum_{i=0}^r h_{r-i}[a/(1-t)]h_i[b/(1-t)] \\
&=a^r\sum_{i=0}^r \frac{(b/a)}{(t;t)_i(t;t)_{r-i}} 
\qquad\qquad\quad \text{(by \eqref{Euler})} \notag \\
&=a^r \frac{H_r(b/a)}{(t;t)_r}. \notag 
\end{align}

Next we turn to Theorem~\ref{thmM}.
Let $Q'_{\la}$ be the modified Hall--Littlewood function
\begin{equation*}
Q'_{\la}[\X;t]=Q_{\la}[\X/(1-t);t].
\end{equation*}
{} From \eqref{Cauchy} it follows that
\begin{equation*}
\sum_{\la}P_{\la}[\X;t]Q'_{\la}[\Y;t]=
\prod_{\substack{x\in\X \\ y\in\Y}}\frac{1}{1-xy}=\sigma_1[\X\Y].
\end{equation*}
Consequently,
\begin{align*}
\sigma_1[\X\Y]P_{\mu}[\X;t]
&=\sum_{\nu}
P_{\nu}[\X;t]P_{\mu}[\X;t]Q'_{\nu}[\Y;t]&& \\
&=\sum_{\la,\nu}
f_{\mu\nu}^{\la}(t)P_{\la}[\X;t]Q'_{\nu}[\Y;t] && \text{(by \eqref{struc})} \\
&=\sum_{\la} P_{\la}[\X;t]Q'_{\la/\mu}[\Y;t] && \text{(by \eqref{skew})}.
\end{align*}
(The above equation also follows by the substitution $\X\to\X/(1-t)$ in an
identity on page 227 of \cite{Macdonald95}).
Summing $\mu$ over the even partitions and replacing $\Y$ by $-\Y$,
we thus find
\begin{equation*}
\sum_{\la} P_{\la}[\X;t]
\sum_{\substack{\mu \\ \mu\text{ even}}} Q'_{\la/\mu}[-\Y;t]=
\sigma_1[-\X\Y]
\sum_{\substack{\mu \\ \mu\text{ even}}}
P_{\mu}[\X;t].
\end{equation*}
Finally note that the sum on the right may be performed by
\eqref{i2}.
Hence we arrive at 
\begin{equation}\label{ehxy}
\sum_{\la} P_{\la}[\X;t]B_{\la}[\Y;t]
=\sigma_1[(1-t)e_2[\X]+h_2[\X]-\X\Y],
\end{equation}
with
\begin{equation*}
B_{\la}[\Y;t]
=\sum_{\substack{\mu \\ \mu\text{ even}}} Q'_{\la/\mu}[-\Y;t].
\end{equation*}
Dispensing with the plethystic notation we may write
\eqref{ehxy} as
\begin{equation*}
\sum_{\la} P_{\la}(x;t)B_{\la}(y;t)
=\prod_{i,j\geq 1}(1-x_iy_j)
\prod_{i\geq 1}\frac{1}{1-x_i^2}
\prod_{i<j}\frac{1-tx_ix_j}{1-x_ix_j},
\end{equation*}
Theorem~\ref{thmM} (with $a\to -a$ and $b\to -b$)
is thus equivalent to the following
closed form expression for $B_{\la}$ in the case of a
two-letter alphabet.
\begin{theorem}\label{thmM2}
Let $\Y=\{a,b\}$. Then
\begin{equation*}
B_{\la}(a,b;t)=
\sum_{\substack{\mu \\ \mu\textup{ even}}} Q'_{\la/\mu}[-\Y;t]
=(-a)^{l(\lao)} h_{\lae}(ab;t)h_{\lao}(b/a;t).
\end{equation*}
\end{theorem}

For example, when $\la=(1^r)$ the only non-vanishing
contribution to the sum over $\mu$ comes from $\mu=0$
(since $Q'_{\la/\mu}$ vanishes if $\mu\not\subset\la$).
Hence
\begin{align*}
B_{(1^r)}(a,b;t)&=Q'_{(1^r)}[-\Y;t] &&\\
&=Q_{(1^r)}[-\Y/(1-t);t] &&\\
&=b_{(1^r)}(t)\:e_r[-\Y/(1-t)] &&
\text{(by \eqref{Qla} and \eqref{eh})} \\
&=(-1)^r (t;t)_r \: h_r[\Y/(1-t)] && \\
&=(-a)^r H_r(b/a)
&& \text{(by \eqref{RShr}),}
\end{align*}
which is in accordance with the right hand side of Theorem~\ref{thmM2}
for $\la=(1^r)$.

\medskip

In much the same way it follows that
\begin{subequations}\label{ms}
\begin{equation}
\sum_{\la} P_{\la}(x;t)C_{\la}(y;t)
=\prod_{i,j\geq 1}\frac{1-tx_iy_j}{1-x_iy_j}
\prod_{i<j}\frac{1-tx_ix_j}{1-x_ix_j},
\end{equation}
with
\begin{equation}\label{msb}
C_{\la}(y;t)
=\sum_{\substack{\mu \\ \mu'\text{ even}}} c_{\mu}(t)Q_{\la/\mu}(y;t),
\end{equation}
\end{subequations}
but as remarked after \eqref{fail}, only for $y=(a)$
does the sum on the right of \eqref{msb} simplify.
Perhaps the best way to understand this case (corresponding to
\eqref{LW1}) is however not through
\eqref{ms} but by adding $a$ to the alphabet $\X$ as explained 
below in the Macdonald polynomial setting.

Consider the identity \eqref{Pc}. For $a=0$ this is \eqref{p349} 
which may be expressed in $\la$-ring notation as 
\begin{equation*}
\sum_{\substack{\la \\ \la'\text{ even}}}
b_{\la}^{\textup{el}}(q,t) P_{\la}[\X;q,t]
=\sigma_1\biggl[\frac{1-t}{1-q} \, e_2[\X]\biggr]=:f[\X].
\end{equation*}
Replacing $\X$ by $\X+\Y$ and using that $e_2[\X+\Y]=e_2[\X]+e_2[\Y]+
e_1[\X\Y]$,
$\sigma_1[\X+\Y]=\sigma_1[\X]\sigma_1[\Y]$
and $P_{\la}[\X+\Y;q,t]=\sum_{\mu} P_{\mu}[\X;q,t]
P_{\la/\mu}[\Y;q,t]$ (for this last result see
\cite[p. 345]{Macdonald95}), this implies
\begin{equation*}
\sum_{\mu} P_{\mu}[\X;q,t]
\sum_{\substack{\la \\ \la'\text{ even}}}
b_{\la}^{\textup{el}}(q,t) P_{\la/\mu}[\Y;q,t]
=f[\X]f[\Y]
\prod_{\substack{x\in\X \\ y\in\Y}}
\frac{(txy;q)_{\infty}}{(xy;q)_{\infty}}.
\end{equation*}
When $\Y$ contains a single letter $a$, so that we have effectively
added $a$ to $\X$, this simplifies 
to
\begin{equation*}
\sum_{\mu} P_{\mu}[\X;q,t]
\sum_{\substack{\la \\ \la'\text{ even} \\
\la-\mu \text{ hor. strip}}}
a^{\abs{\la-\mu}}
b_{\la}^{\textup{el}}(q,t)\psi_{\la/\mu}(q,t)
=f[\X]
\prod_{x\in\X}
\frac{(atx;q)_{\infty}}{(ax;q)_{\infty}}.
\end{equation*}
To get the expression on the left we have used 
\cite[p. 346]{Macdonald95}.
The partition $\la$ in the sum on the left 
is fixed by $\mu$ as $\la'_i=2\lceil \mu'_i/2\rceil$.
Assuming such $\mu$, we get
\begin{equation*} 
\sum_{\mu}a^{l((\la')_{\textup{o}})} 
b_{\la}^{\textup{el}}(q,t)\psi_{\la/\mu}(q,t)
P_{\mu}[\X;q,t]
=f[\X] \prod_{x\in\X}
\frac{(atx;q)_{\infty}}{(ax;q)_{\infty}}.
\end{equation*}
But putting together two combinatorial identities 
on pages 350 and 351 of \cite{Macdonald95} yields
\begin{equation*} 
b_{\la}^{\textup{el}}(q,t)\psi_{\la/\mu}(q,t)
=b_{\mu}^{\textup{el}}(q,t)
\end{equation*}
so that \eqref{Pc} follows.

\subsection*{Acknowledgements}
I thank the anonymous referee for helpful suggestions,
leading to the material of Section~\ref{sec4}.

\bibliographystyle{amsplain}

\end{document}